\pdfoutput=1
\documentclass[12pt,a4paper]{article}
\usepackage[T1]{fontenc}
\usepackage[left=1in, right=1in]{geometry}
\usepackage{amsthm}
\usepackage{amsmath}
\usepackage{amsxtra}
\usepackage{amstext}
\usepackage{amssymb}
\usepackage[colorlinks, citecolor=blue, linkcolor=blue]{hyperref}
\theoremstyle{definition}
\newtheorem{defn}{Definition}[section]
\newtheorem{rem}{Remark}[section]
\newtheorem{eg}{Example}[section]
\theoremstyle{plain}
\newtheorem{theo}{Theorem}[section]
\newtheorem{lem}{Lemma}[section]
\newtheorem{cor}[lem]{Corollary}
\numberwithin{equation}{section}
\begin{document}
		\title{\textbf{Pointwise Bi-Slant Submanifolds in Locally Conformal Kähler Manifolds Immersed as Warped Products}}
		\author{\textbf{Umar Mohd Khan \& Viqar Azam Khan}\\Department of Mathematics\\Aligarh Muslim University\\Aligarh-202002, India\\Email: umar.007.morpheus@gmail.com, viqarster@gmail.com}
		\date{}
		\maketitle
		\begin{abstract}
			We study immersions of pointwise bi-slant submanifolds of locally conformal Kähler manifolds as warped products. In particular, we establish characterisation theorem for a pointwise bi-slant submanifold of a locally conformal Kähler manifold to be immersed as a warped product and show that a necessary condition is that the Lee vector field $B$ is orthogonal to the second factor and the warping function $\lambda$ satisfies $\text{grad}(\ln\lambda)=\frac{1}{2}B^T$, where $B^T$ denotes the tangential part of the Lee vector field. We also extend Chen's inequality for the squared length of the second fundamental form to our case and study the corresponding equality case.
		\end{abstract}
	\section{Introduction}
	Vaisman introduced locally conformal Kähler (lcK) manifolds as a generalisation of Kähler manifolds \cite{Vaisman-lcK1,Vaisman-lcK2,Vaisman-lcK3,Vaisman-lcK4,Vaisman-lcK5,Vaisman-lcK6,Vaisman-lcK7}. An lcK manifold is a Hermitian manifold that can be written as the union of Kähler manifolds such that the lcK metric is locally conformal to these Kähler metrics. LcK manifolds are characterised by the existence of a globally defined closed 1-form $\omega$, called the \textit{Lee form}, such that the fundamental 2-form of the lcK metric satisfies $d\Omega=\Omega\wedge\omega$. The Lee form and its associated Lee Vector field play an important part in the geometry of lcK manifolds.\par
	From an extrinsic geometric standpoint, holomorphic and totally real submanifolds are important objects of study in the setting of almost Hermitian manifolds. Bejancu \cite{Bejancu-CR-I,Bejancu-CR-II} defined CR submanifolds as a generalisation of holomorphic and totally real submanifolds which were further studied by Chen \cite{Chen-CR-I,Chen-CR-II}. Later, Chen \cite{Chen-SlantImm,Chen-SlantSubmanifolds} extended the class of holomorphic and totally real submanifolds by introducing the notion of slant submanifolds. The concept was further generalised to pointwise slant submanifolds \cite{Chen-PointwiseSlant} by the same author. The study of CR submanifolds and slant submanifolds was later generalised by several authors to semi-slant submanifolds, hemi-slant submanifolds(also called pseudo-slant submanifolds) and bi-slant submanifolds, in various ambient manifolds.\par
	Semi-slant submanifolds in almost Hermitian manifolds were studied by Papaghiuc \cite{Papaghiuc-SemiSlant}. Cabrerizo et al. \cite{Cabrerizo-SemiSlant-Sasakian,Cabrerizo-Slant-Sasakian} studied semi-slant submanifolds in Sasakian manifolds. Slant and semi-slant submanifolds in almost product Riemannian manifolds were studied in \cite{Atceken-Slant-RiemProd,Li_Liu-SemiSlant-LocProd,Sahin-Slant-AlmProdRiem}. Hemi-slant submanifolds were also studied in nearly Kenmotsu manifolds \cite{Atceken-HemiSlant-Kenmotsu}, LCS-manifolds \cite{Atceken-Slant_HemiSlant-LCS} and locally product Riemannian manifolds \cite{Tastan_Ozemdir-HemiSlant-LocProdRiem}.\par
	Bishop and O'Neill \cite{WarpedProd-Def} while studying examples of manifolds with negative sectional curvature, defined warped product manifolds by homothetically warping the product metric on a product manifold. Warped products are a natural generalisation of Riemannian products and they have found extensive applications in relativity. Most notably the Schwarzschild metric describing the gravitational field outside a spherical mass under certain assumptions and the Robertsen Walker metric (FLRW metric) are examples of warped product metrics. A natural example of warped product manifolds are surfaces of revolution. Hiepko \cite{Hiepko-WarpedProdCharacterisation} gave a characterisation for a Riemannian manifold to be the warped product of its submanifolds, generalising the deRham decomposition theorem for product manifolds. Later on Nölker \cite{Nolker-WarpedExtrinsic} and Chen \cite{Chen-TwistedProd,Chen-WarpedProd1,Chen-WarpedProd2} initiated the study of extrinsic geometry of warped product manifolds.\par
	Chen \cite{Chen-WarpedProdCR-CRWarpedProd-I,Chen-WarpedProdCR-CRWarpedProd-II} initiated the study of CR submanifolds immersed as warped products in Kähler manifolds. He proved that given any holomorphic $M_T$ and totally real submanifold $M_\perp$ of a Kähler manifold, every warped product of the form $M_T\times_\lambda M_\perp$ in a Kähler manifold satisfies the inequality
	\begin{equation}\label{Intro:Chen'sIneq}
		||h||^2\geq 2n_2||\text{grad}(\ln\lambda)||^2
	\end{equation} 
	where $\lambda$ is the warping function, $n_2$ is the dimension of $M_\perp$, $||h||^2$ is the squared norm of the second fundamental form and $\text{grad}(\ln\lambda)$ is the gradient of $\ln\lambda$. Bonanzinga and Matsumoto \cite{Bonanzinga_Matsumoto-WarpedProdCR-lcK,Matsumoto_Bonanzinga-DoublyWarpedProdCR-lcKSpaceForm,Matsumoto_Bonanzinga-DoublyWarpedProdCR-lcKSpaceForm-II} continued the study in the setting of lcK manifolds. Nargis Jamal et al. \cite{NargisJamal_KAK_VAK-GenericWarped-lcK} studied Generic warped products in lcK manifolds. Further studies of semi-slant and hemi-slant submanifolds of lcK manifolds were carried out in \cite{Alghamdi_Uddin-SemiSlantWarped-lcK,Tastan_Sibel-HemiSlant-lcK,Tastan_Tripathi-SemiSlant-lcK}. Generic submanifolds, CR-submanifolds and pointwise semi-slant submanifolds immersed as warped products in lcK manifolds were studied by \cite{NargisJamal_KAK_VAK-GenericWarped-lcK,Alghamdi_Uddin-SemiSlantWarped-lcK}.\par
	We continue the study by considering pointwise bi-slant submanifolds in an lcK manifold. In particular we give characterisation theorems and establish Chen's inequality for the squared norm of the second fundamental form of pointwise bi-slant submanifolds immersed as warped products in an lcK manifold.
	\section{Preliminaries}
	\allowdisplaybreaks
	\begin{defn}\label{def:lcK manifold}
		A Hermitian Manifold $(\widetilde{M}^{2n},J,g )$ is said to be a \textit{locally conformal Kähler} (lcK) manifold if there exists an open cover $\{U_i\}_{i\in I}$ of $\widetilde{M}^{2n}$ and a family $\{f_i\}_{i\in I}$ of $C^\infty$ functions $f_i:U_i\to\mathbb{R}$ such that for each $i\in I$, the metric 
		\begin{equation}\label{eq:g_i defn}
			g_i=e^{-f_i}g|_{U_i} 
		\end{equation} 
		on $U_i$ is a Kähler metric.
	\end{defn}
	Given an lcK manifold $(\widetilde{M}^{2n},J,g )$, let $U,V$ denote smooth sections of $T\widetilde{M}^{2n}$, then the local 1-forms $df_i$ glue up to a globally defined closed 1-form $\omega$, called the \textit{Lee form}, and it satisfies the following equation
	\begin{equation}\label{eq:d(Fund 2-form)}
		d\Omega=\Omega\wedge\omega
	\end{equation} 
	where $\Omega(U,V)=g (JU,V)$ is the fundamental 2-form associated to $(J,g )$.\par
	Denote by $B$ the vector field equivalent to $\omega$ with respect to $g$, i.e. $\omega(U)=g(B,U)$. $B$ is called the \textit{Lee vector field}.\par
	Let $\overline{\nabla}$ denote the Levi-Civita connection of $(\widetilde{M}^{2n},g)$ and $\widetilde{\nabla}_i$ denote the Levi-Civita connection of the local metrics $g_i$ for all $i\in I$. Then $\widetilde{\nabla}_i$ glue up to a globally defined torsion-free linear connection $\widetilde{\nabla}$ on $\widetilde{M}^{2n}$ given by
	\begin{equation}\label{eq:Weyl Conn defn}
		\widetilde{\nabla}_UV=\overline{\nabla}_UV-\frac{1}{2}\left\lbrace\omega(U)V+\omega(V)U-g(U,V )B \right\rbrace  
	\end{equation}
	where $U,V\in T\widetilde{M}^{2n}$ and satisfying 
	\begin{equation}\label{eq:Weyl Conn and g}
		\widetilde{\nabla}g=\omega\otimes g
	\end{equation}
	$\widetilde{\nabla}$ is called the \textit{Weyl connection} of the lcK manifold $(\widetilde{M}^{2n},J,g )$. As  $g_i$ are Kähler metrics, the almost complex structure $J$ is parallel with respect to the Weyl connection, i.e. $\widetilde{\nabla}J=0$. This gives 
	\begin{equation}\label{eq:Riemm Conn and J}
		\overline{\nabla}_UJV=J\overline{\nabla}_UV+\frac{1}{2}\left\lbrace \Theta(V)U-\omega(V)JU-g(U,V)A+\Omega(U,V)B\right\rbrace 
	\end{equation}
	Now as $\omega$ is a closed form on $\widetilde{M}^{2n}$, we have 
	\begin{equation}\label{eq:Derivative Lee Form Riemm Conn}
		(\overline{\nabla}_U\omega)V=(\overline{\nabla}_V\omega)U
	\end{equation}
	Let $M^m$ be a Riemannian manifold isometrically immersed in an lcK manifold $(\widetilde{M}^{2n},J,g )$. Let $U,V,W$ denote smooth sections of $TM^m$ and $\xi,\eta$ denote smooth sections of $T^\perp M^m$. \par
	The Gauss and Weingarten formulae with respect to the Riemannian connection of $\widetilde{M}^{2n}$ are given as 
	\begin{align}
		\overline{\nabla}_UV&=\nabla_UV+h(U,V)\label{eq:Gauss - Riemm Conn}\\
		\overline{\nabla}_U\xi&=-\mathfrak{A}_\xi U+\nabla^\perp_U \xi\label{eq:Weingarten - Riemm Conn}
	\end{align}
	where $h$ is the second fundamental form, $\mathfrak{A}$ is the shape operator and $\nabla, \nabla^\perp$ are respectively the induced connections in the tangent bundle and the normal bundle of $M^m$ with respect to $\overline{\nabla}$.\par  
	The Gauss and Weingarten formulae with respect to the Weyl connection of $\widetilde{M}^{2n}$ are given as 
	\begin{align}
		\widetilde{\nabla}_UV&=\hat{\nabla}_UV+\widetilde{h}(U,V)\label{eq:Gauss - Weyl Conn}\\
		\widetilde{\nabla}_U\xi&=-\widetilde{\mathfrak{A}}_\xi U+\widetilde{\nabla}^\perp_U \xi\label{eq:Weingarten - Weyl Conn}
	\end{align}
	where $\widetilde{h}$ is the second fundamental form, $\widetilde{\mathfrak{A}}$ is the shape operator and $\hat{\nabla}, \widetilde{\nabla}^\perp$ are respectively the induced connections in the tangent bundle and the normal bundle of $M^m$ with respect to $\widetilde{\nabla}$.\par
	Let $H$ denote the trace of $h$, then $H$ is called the mean curvature vector of $M^m$ in $(\widetilde{M}^{2n},J,g )$ and is a smooth section of $T^\perp M^m$. We say $M^m$ is a totally umbilic submanifold of $(\widetilde{M}^{2n},J,g )$, if $h(U,V)=g(U,V)H$. We say $M^m$ is a totally geodesic submanifold of $(\widetilde{M}^{2n},J,g )$, if $h(U,V)=0$.\par
	Let $B^T, B^N$ denote the tangential and normal components of the Lee vector field $B$.\par
	From \eqref{eq:Weyl Conn defn}, we have the following relations 
	\begin{align}
		\hat{\nabla}_UV&=\nabla_UV-\frac{1}{2}\left\lbrace\omega(U)V+\omega(V)U-g(U,V)B^T  \right\rbrace \label{eq:Relation induced conn M - Weyl and Riemm}\\
		\widetilde{h}(U,V)&=h(U,V)+\frac{1}{2}g(U,V)B^N \label{eq:Relation 2nd fund form - Weyl and Riemm}\\
		\widetilde{\mathfrak{A}}_\xi U&=\mathfrak{A}_\xi U+\frac{1}{2}\omega(\xi)U \label{eq:Relation shape operator - Weyl and Riemm}\\
		\widetilde{\nabla}^\perp_U\xi&=\nabla^\perp_U\xi-\frac{1}{2}\omega(U)\xi \label{eq:Relation induced normal conn - Weyl and Riemm}
	\end{align}
	Now define
	\begin{align}\label{eq:Definition P,F,t,f}
		JU&=PU+FU&J\xi=t\xi+f\xi
	\end{align}
	where $PU, t\xi$ and $FU, f\xi$ are respectively the tangential and normal parts. Then, we have
	\begin{equation}\label{eq:Identities P,F,t,f}
		\begin{aligned}
			P^2+tF&=-I&\hspace{1cm}f^2+Ft&=-I\\
			FP+fF&=0&\hspace{1cm}tf+Pt&=0
		\end{aligned}
	\end{equation}
	Define the covariant differentiation of $P$, $F$, $t$ and $f$ with respect to the Levi-Civita connection of $\widetilde{M}^{2n}$ as 
	\begin{equation}\label{eq:Definition Covar Diff Riem Conn P,F,t,f}
		\begin{aligned}
			(\overline{\nabla}_UP)V& =\nabla_UPV-P\nabla_UV\\
			(\overline{\nabla}_UF)V& =\nabla^\perp_UFV-F\nabla_UV\\
			(\overline{\nabla}_Ut)\xi& =\nabla_ut\xi-t(\nabla^\perp_U\xi)\\
			(\overline{\nabla}_Uf)\xi& =\nabla^\perp_Uf\xi-f(\nabla^\perp_U\xi)
		\end{aligned}
	\end{equation}
	Then as $\widetilde{\nabla}J=0$, using \eqref{eq:Relation induced conn M - Weyl and Riemm}, \eqref{eq:Relation 2nd fund form - Weyl and Riemm}, \eqref{eq:Relation shape operator - Weyl and Riemm}, \eqref{eq:Relation induced normal conn - Weyl and Riemm} we have 
	\begin{equation}\label{eq:Relation Covar Diff Riem Conn P,F,t,f}
		\begin{aligned}
			(\overline{\nabla}_UP)V&=\mathfrak{A}_{FV}U+th(U,V)+\frac{1}{2}\left\lbrace \Theta(V)U-\omega(V)PU+g(PU,V)B^T-g(U,V)A^T\right\rbrace\\
			(\overline{\nabla}_UF)V&=fh(U,V)-h(U,PV)+\frac{1}{2}\left\lbrace g(PU,V)B^N-g(U,V)A^N-\omega(V)FU\right\rbrace\\
			(\overline{\nabla}_Ut)\xi&=\mathfrak{A}_{f\xi}U-P\mathfrak{A}_\xi U+\frac{1}{2}\left\lbrace g(FU,\xi)B^T-\omega(\xi)PU+\Theta(\xi)U \right\rbrace\\
			(\overline{\nabla}_Uf)\xi&=-h(U,t\xi)-F\mathfrak{A}_\xi U+\frac{1}{2}\left\lbrace g(FU,\xi)B^N-\omega(\xi)FU\right\rbrace 
		\end{aligned}
	\end{equation}
	Bishop and O'Neill \cite{WarpedProd-Def} defined warped product as
	\begin{defn}\label{defn:Warped Product Manifolds}
		Let $(M_1^{n_1},g_1)$ and $(M_2^{n_2},g_2)$ be Riemmanian manifolds and let 
		\[\pi_1:M_1\times M_2\to M_1\text{ and }\pi_2:M_1\times M_2\to M_2 \]
		be the canonical projections. Let \mbox{$\lambda:M_1\to(0,\infty)$} be a smooth function. Then the warped product manifold \mbox{$(M,g)=M_1\times \,_\lambda M_2$} is defined as the manifold $M_1\times M_2$ equipped with the Riemannian metric
		\begin{equation}\label{eq:Warped Product Metric}
			g=\pi_1^\star g_1+\lambda^2\pi_2^\star g_2
		\end{equation}
	\end{defn}
	\noindent Warped product manifolds are a generalization of the usual product of two Riemannian manifolds. In fact we have the following characterisation theorem
	\begin{theo}[\cite{Hiepko-WarpedProdCharacterisation}]\label{th:Characterisation - Warped Product}
		Let $(M^m,g)$ be a connected Riemannian manifold equipped with orthogonal, complementary, involutive distributions $\mathcal{D}_1$ and $\mathcal{D}_2$. Further let the leaves of $\mathcal{D}_1$ be totally geodesic and the leaves of $\mathcal{D}_2$ be extrinsic spheres in $M^m$, where by extrinsic spheres we mean totally umbilic submanifolds such that the mean curvature vector is parallel in the normal bundle. Then $(M^m,g)$ is locally a warped product \mbox{$(M,g)=M_1\times \,_\lambda M_2$}, where $M_1$ and $M_2$ respectively denote the leaves of $\mathcal{D}_1$ and $\mathcal{D}_2$ and \mbox{$\lambda:M_1\to(0,\infty)$} is a smooth function such that $\text{grad}(\ln\lambda)$ is the mean curvature vector of $M_2$ in $M$.\par
		Further, if $(M^m,g)$ is simply connected and complete, then $(M^m,g)$ is globally a warped product.
	\end{theo}
	\noindent For $(M_1^{n_1},g_1)$, $(M_2^{n_2},g_2)$ and $(M,g)$ denote respectively the Levi-Civita connections by $\nabla^1$, $\nabla^2$ and $\nabla$. Given any smooth function $\lambda:M_1\to\mathbb{R}$, let $\text{grad}(\lambda)$ denote the lift of the gradient vector field of $\lambda$ to $(M,g)$.
	\begin{theo}[\cite{Hiepko-WarpedProdCharacterisation}]\label{th:Expressions Connection Warped Product}
		Given a warped product manifold \mbox{$(M,g)= M_1\times \,_\lambda M_2$} of Riemmanian manifolds $(M_1^{n_1},g_1)$ and $(M_2^{n_2},g_2)$, we have for all $X,Y\in \mathcal{L}(M_1)$ and $Z,W\in \mathcal{L}(M_2)$,
		\begin{align}
			\nabla_XY&=\nabla^1_XY\label{eq:Connection Warped Product on M_1 x M_1}\\
			\nabla_XZ&=\nabla_ZX=X(\ln\lambda)Z\label{eq:Connection Warped Product on M_1 x M_2}\\
			\nabla_ZW&=\nabla^2_ZW-g(Z,W)\text{grad}(\ln\lambda)\label{eq:Connection Warped Product on M_2 x M_2}
		\end{align}
	\end{theo}
	\noindent It follows from Lemma \ref{th:Expressions Connection Warped Product} that $\mathcal{H}=-\text{grad}(\ln\lambda)$ is the mean curvature vector of $M_2$ in $M$.
	
	\section{Pointwise Bi-Slant Submanifolds of lcK manifolds}
	Let $M^m$ be a Riemannian manifold isometrically immersed in an lcK manifold $(\widetilde{M}^{2n},J,g )$. $M^m$ is said to be a \textit{pointwise bi-slant submanifold} if it admits two orthogonal complementary distributions $\mathcal{D}^{\theta_1}$ and $\mathcal{D}^{\theta_2}$, such that $\mathcal{D}^{\theta_1}$ and $\mathcal{D}^{\theta_2}$ are pointwise slant with slant angles $\theta_1,\theta_2\in\left(0,\frac{\pi}{2}\right)$ and $\theta_1\neq\theta_2$, i.e. $P^2X=-\cos^2\theta_1 X$, for every smooth vector field $X\in\mathcal{D}^{\theta_1}$ and $P^2Z=-\cos^2\theta_2 Z$, for every smooth vector field $Z\in\mathcal{D}^{\theta_2}$.\par
	The tangent bundle and the normal bundle of a pointwise bi-slant submanifold admits an orthogonal decomposition as
	\begin{align}\label{eq:Bi Slant Tangent, Normal Bundle Ortho decomp }
		TM^m&=\mathcal{D}^{\theta_1}\oplus\mathcal{D}^{\theta_2}&T^\perp M^m&=F\mathcal{D}^{\theta_1}\oplus F\mathcal{D}^{\theta_2}+\mu
	\end{align}
	where $\mu$ is the orthogonal complementary distribution of $F\mathcal{D}^{\theta_1}\oplus F\mathcal{D}^{\theta_2}$ in $T^\perp M^m$ and is an invariant subbundle of $T^\perp M^m$ with respect to $J$. It is easy to observe that for $i=1,2$,
	\begin{equation}\label{eq:Bi-Slant - Image of D,D' under P,F,t}
		\begin{aligned}
			P\mathcal{D}^{\theta_i}&=\mathcal{D}^{\theta_i}&t(F\mathcal{D}^{\theta_i})&=\mathcal{D}^{\theta_i}&\hspace{1cm}t(\mu)&=\{0\}\\
			f(F\mathcal{D}^{\theta_i})&= F\mathcal{D}^{\theta_i}&\hspace{1cm}f(\mu)&=\mu
		\end{aligned}
	\end{equation}
	Let $M^m$ be a pointwise bi-slant manifold isometrically immersed in an lcK manifold $(\widetilde{M}^{2n},J,g )$ such that the distributions $\mathcal{D}^{\theta_1}, \mathcal{D}^{\theta_2}$ are both involutive. Let $M_{\theta_1}^{2n_1}$ and $M_{\theta_2}^{2n_2}$ respectively denote the leaves of $\mathcal{D}^{\theta_1}$ and $\mathcal{D}^{\theta_2}$, where $2n_1=\dim_\mathbb{R}\mathcal{D}^{\theta_1}$ and $2n_2=\dim_\mathbb{R}\mathcal{D}^{\theta_2}$. We say $M^m$ is a 
	\begin{itemize}
		\item \textit{mixed totally geodesic pointwise bi-slant submanifold} if $h(\mathcal{D}^{\theta_1},\mathcal{D}^{\theta_2})=\{0\}$.
		\item \textit{pointwise bi-slant product submanifold} if $M^m$ can be expressed locally as \mbox{$M_{\theta_1}\times M_{\theta_2}$}.
		\item \textit{pointwise bi-slant warped product submanifold} if $M^m$ can be expressed locally as \mbox{$M_{\theta_1}\times\,_\lambda M_{\theta_2}$} for some smooth function $\lambda:M_{\theta_1}\to(0,\infty)$.
	\end{itemize}
	Let $X,Y$ be smooth vector fields in $\mathcal{D}^{\theta_1}$ and $Z,W$ be smooth vector fields in $\mathcal{D}^{\theta_2}$. Then we have,
	\begin{theo}\label{th:Bi-Slant D_theta1}
		Let $M^m$ be a pointwise bi-slant submanifold of an lcK manifold $\widetilde{M}^{2n}$. 
		\begin{itemize}
			\item the slant distribution $\mathcal{D}^{\theta_1}$ is involutive if and only if 
			\begin{equation}\label{eq:bi-slant D_theta1 involutive}
				\begin{aligned}
					g(\mathfrak{A}_{FPX}Z-\mathfrak{A}_{FX}PZ,Y)&-g(\mathfrak{A}_{FPY}Z-\mathfrak{A}_{FY}PZ,X)\\
					&=g(\nabla^\perp_XFY,FZ)-g(\nabla^\perp_YFX,FZ)
				\end{aligned}
			\end{equation}
			\item the leaves of the slant distribution $\mathcal{D}^{\theta_1}$ are totally geodesic in $M^m$ if and only if 
			\begin{equation}\label{eq:bi-slant D_theta1 totally geodesic}
				\omega(\mathcal{D}^{\theta_2})=\{0\}\text{ and }g(\mathfrak{A}_{FPX}Z-\mathfrak{A}_{FX}PZ,Y)+g(\nabla^\perp_YFX,FZ)=0
			\end{equation}
			\item the leaves of the slant distribution $\mathcal{D}^{\theta_1}$ are totally umbilic in $M^m$ if and only if 
			\begin{equation}\label{eq:bi-slant D_theta1 totally umbilic}
				\begin{aligned}
					g(\mathfrak{A}_{FPX}Z-\mathfrak{A}_{FX}PZ,Y)&+g(\nabla^\perp_YFX,FZ)\\
					&=\sin^2\theta_1\left(\frac{1}{2}\omega(Z)+g(\mathcal{H},Z)\right)g(X,Y)
				\end{aligned}
			\end{equation}
			for some smooth vector field $\mathcal{H}\in\mathcal{D}^{\theta_2}$.
		\end{itemize}
	\end{theo}
	\begin{proof}
		From \eqref{eq:Riemm Conn and J} and \eqref{eq:Identities P,F,t,f}, we have 
		\begin{align*}
			g(\nabla_XY,Z)&=g(\overline{\nabla}_XPY+\overline{\nabla}_XFY-\frac{1}{2}g(X,Y)JB-\frac{1}{2}g(PX,Y)B,JZ)\\
			&=-g(J\overline{\nabla}_XPY,Z)-g(\mathfrak{A}_{FY}X,PZ)+g(\nabla^\perp_XFY,FZ)\\
			&\hspace{0.4cm}-\frac{1}{2}g(X,Y)g(B,Z)-\frac{1}{2}g(PX,Y)g(B,JZ)\\
			&=-g(\overline{\nabla}_XJPY,Z)+\frac{1}{2}g(PX,PY)g(B,Z)-g(\mathfrak{A}_{FY}X,PZ)\\
			&\hspace{0.4cm}-\frac{1}{2}g(X,Y)g(B,Z)+g(\nabla^\perp_XFY,FZ)\\
			&=\cos^2\theta_1 g(\nabla_XY,Z)+g(\mathfrak{A}_{FPY}X,Z)+\frac{1}{2}\sin^2\theta_1 g(X,Y)g(B,Z)\\
			&\hspace{0.4cm}-g(\mathfrak{A}_{FY}X,PZ)+g(\nabla^\perp_XFY,FZ)
		\end{align*}
		i.e. we have
		\[\sin^2\theta_1\left(g(\nabla_XY,Z)+\frac{1}{2}g(X,Y)g(B,Z)\right)=g(\mathfrak{A}_{FPY}Z-\mathfrak{A}_{FY}PZ,X)+g(\nabla^\perp_XFY,FZ)\]
		Hence, the result follows.
	\end{proof}	
	Similarly, we have
	\begin{theo}\label{th:Bi-Slant D_theta2}
		Let $M^m$ be a pointwise bi-slant submanifold of an lcK manifold $\widetilde{M}^{2n}$. Then 
		\begin{itemize}
			\item the slant distribution $\mathcal{D}^{\theta_2}$ is involutive if and only if 
			\begin{equation}\label{eq:bi-slant D_theta2 involutive}
				\begin{aligned}
					g(\mathfrak{A}_{FPZ}X-\mathfrak{A}_{FZ}PX,W)&-g(\mathfrak{A}_{FPW}X-\mathfrak{A}_{FW}PX,Z)\\
					&=g(\nabla^\perp_ZFW,FX)-g(\nabla^\perp_WFZ,FX)
				\end{aligned}
			\end{equation}
			\item the leaves of the slant distribution $\mathcal{D}^{\theta_2}$ are totally geodesic in $M^m$ if and only if 
			\begin{equation}\label{eq:bi-slant D_theta2 totally geodesic}
				\omega(\mathcal{D}^{\theta_1})=\{0\}\text{ and }g(\mathfrak{A}_{FPZ}X-\mathfrak{A}_{FZ}PX,W)+g(\nabla^\perp_WFZ,FX)=0
			\end{equation}
			\item the leaves of the slant distribution $\mathcal{D}^{\theta_2}$ are totally umbilic in $M^m$ if and only if 
			\begin{equation}\label{eq:bi-slant D_theta2 totally umbilic}
				\begin{aligned}
					g(\mathfrak{A}_{FPZ}X-\mathfrak{A}_{FZ}PX,W)&+g(\nabla^\perp_WFZ,FX)\\
					&=\sin^2\theta_2\left(\frac{1}{2}\omega(X)+g(\mathcal{H},X)\right)g(Z,W)
				\end{aligned}
			\end{equation}
			for some smooth vector field $\mathcal{H}\in\mathcal{D}^{\theta_1}$.
		\end{itemize}
	\end{theo}
	\noindent \textbf{Notations:} Let $\mathcal{D}^{\theta_1}$ and $\mathcal{D}^{\theta_2}$ be the slant distributions on a pointwise bi-slant submanifold $M^m$ of an lcK manifold $\widetilde{M}^{2n}$ such that both distributions are involutive and let $M_{\theta_1}$ and $M_{\theta_2}$ respectively denote the leaves of the distributions $\mathcal{D}^{\theta_1}$ and $\mathcal{D}^{\theta_2}$ respectively. Then $\mathcal{D}^{\theta_1}(p,q)=T_{(p,q)}(M_{\theta_1}\times\{q\})$ and $\mathcal{D}^{\theta_2}(p,q)=T_{(p,q)}(\{p\}\times M_{\theta_2})$. Let $\mathcal{L}(M_{\theta_1})$ and $\mathcal{L}(M_{\theta_2})$ respectively denote the set of lifts of vector fields from $M_{\theta_1}$ and $M_{\theta_2}$ to $M$. Then $X\in\mathcal{L}(M_{\theta_1})$ if and only if $X|_{\{p\}\times M_{\theta_2}}$ is constant for every $p\in M_{\theta_1}$. Similarly, $Z\in\mathcal{L}(M_{\theta_2})$ if and only if $Z|_{M_{\theta_1}\times\{q\}}$ is constant for every $q\in M_{\theta_2}$. Also, if \mbox{$\pi_{\theta_1}:M_{\theta_1}\times M_{\theta_2}\to M_{\theta_1}$} and \mbox{$\pi_{\theta_2}:M_{\theta_1}\times M_{\theta_2}\to M_{\theta_2}$} are the canonical projections, we have $d\pi_{\theta_1}(\mathcal{L}(M_{\theta_1}))=TM_{\theta_1}$ and $d\pi_{\theta_2}(\mathcal{L}(M_{\theta_2}))=TM_{\theta_2}$. It is clear that a general vector field in $\mathcal{D}^{\theta_1}$ (respectively $\mathcal{D}^{\theta_2}$) need not be in $\mathcal{L}(M_{\theta_1})$ (respectively $\mathcal{L}(M_{\theta_2})$).\par
	From here on we use $X,Y$ to denote smooth vector fields in $\mathcal{L}(M_{\theta_1})$ and $Z,W$ to denote smooth vector fields in $\mathcal{L}(M_{\theta_2})$
	\section{Some Lemmas}
	We give the following lemmas which will be used to prove our main results.
	\begin{lem}\label{lem:Identities For Bi-Slant Warped Product Submanifolds}
		Given a pointwise bi-slant warped product submanifold \mbox{$M= M_{\theta_1}\times \,_{\lambda}M_{\theta_2}$} in an lcK manifold $(\widetilde{M}^{2n},J,g )$, we have for all $X,Y\in \mathcal{L}(M_{\theta_1})$ and $Z,W\in \mathcal{L}(M_{\theta_2})$, 
		\begin{align}
			g(h(X,Z),FY)&=g(h(Y,Z),FX)\label{eq:g(h(D_1,D_2),FD_1) Symmetric Bi-Slant Warped Product Submanifolds}\\
			g(h(X,Z),FW)&=g(h(X,W),FZ)\label{eq:g(h(D_1,D_2),FD_2) Symmetric Bi-Slant Warped Product Submanifolds}\\
			g(h(X,Y),FZ)&=g(h(X,Z),FY)-\frac{1}{2}g(X,Y)g(B,FZ)\label{eq:Identity g(h(D_1,D_1),FD_2) & g(h(D_1,D_2),FD_1) Bi-Slant Warped Product Submanifolds}\\
			g(h(Z,W),FX)&=g(h(X,Z),FW)-\frac{1}{2}g(Z,W)g(B,FX)\label{eq:Identity g(h(D_2,D_2),FD_1) & g(h(D_1,D_2),FD_2) Bi-Slant Warped Product Submanifolds}\\
			X(\ln\lambda)&=\frac{1}{2}g(B,X)\label{eq:Identity g(B,D_1) Bi-Slant Warped Product Submanifolds}\\
			g(B,Z)&=0\label{eq:Identity g(B,D_2) Bi-Slant Warped Product Submanifolds}
		\end{align}
	\end{lem}
	\begin{proof}
		For all $X,Y\in \mathcal{L}(M_{\theta_1})$ and $Z,W\in \mathcal{L}(M_{\theta_2})$, we have using \eqref{eq:Riemm Conn and J} and \eqref{eq:Connection Warped Product on M_1 x M_2},
		\begin{align*}
			g(h(X,Z),FW)&=g(\overline{\nabla}_XZ,JW-PW)\\
			&=-g(J\overline{\nabla}_XZ,W)-g(\nabla_XZ,PW)\\
			&=-g(\overline{\nabla}_XJZ,W)-X(\ln\lambda)g(Z,PW)\\
			&=-g(\overline{\nabla}_XPZ,W)-g(\overline{\nabla}_XFZ,W)-X(\ln\lambda)g(Z,PW)\\
			&=-X(\ln\lambda)g(PZ,W)+g(\mathfrak{A}_{FZ}X,W)-X(\ln\lambda)g(Z,PW)\\
			&=g(h(X,W),FZ)
			\intertext{which gives \eqref{eq:g(h(D_1,D_2),FD_2) Symmetric Bi-Slant Warped Product Submanifolds}. Repeating the above calculation, we have}
			g(h(X,Z),FW)&=g(\overline{\nabla}_ZX,JW-PW)\\
			&=-g(J\overline{\nabla}_ZX,W)-g(\nabla_ZX,PW)\\
			&=-g(\overline{\nabla}_ZJX,W)-\frac{1}{2}g(JB,X)g(Z,W)-\frac{1}{2}g(B,X)g(JZ,W)\\
			&\hspace{0.4cm}-X(\ln\lambda)g(Z,PW)\\
			&=-g(\overline{\nabla}_ZPX,W)-g(\overline{\nabla}_ZFX,W)-\frac{1}{2}g(JB,X)g(Z,W)\\
			&\hspace{0.4cm}-\frac{1}{2}g(B,X)g(JZ,W)-X(\ln\lambda)g(Z,PW)\\
			&=-PX(\ln\lambda)g(Z,W)+g(\mathfrak{A}_{FX}Z,W)+\frac{1}{2}g(B,JX)g(Z,W)\\
			&\hspace{0.4cm}-\frac{1}{2}g(B,X)g(PZ,W)-X(\ln\lambda)g(Z,PW)
			\intertext{Using \eqref{eq:g(h(D_1,D_2),FD_2) Symmetric Bi-Slant Warped Product Submanifolds} and comparing symmetric and skew symmetric terms in $Z$ and $W$ we have,}
			X(\ln\lambda)&=\frac{1}{2}g(B,X)
			\intertext{which gives \eqref{eq:Identity g(B,D_1) Bi-Slant Warped Product Submanifolds} and }
			g(h(X,Z),FW)&=g(h(Z,W),FX)-PX(\ln\lambda)g(Z,W)+\frac{1}{2}g(B,PX+FX)g(Z,W)
			\intertext{which on substituting from \eqref{eq:Identity g(B,D_1) Bi-Slant Warped Product Submanifolds} gives \eqref{eq:Identity g(h(D_2,D_2),FD_1) & g(h(D_1,D_2),FD_2) Bi-Slant Warped Product Submanifolds}. Similarly,}
			g(h(X,Z),FY)&=g(\overline{\nabla}_ZX,JY-PY)\\
			&=-g(J\overline{\nabla}_ZX,Y)-g(\nabla_ZX,PY)\\
			&=-g(\overline{\nabla}_ZJX,Y)\\
			&=-g(\overline{\nabla}_ZPX,Y)-g(\overline{\nabla}_ZFX,Y)\\
			&=g(\mathfrak{A}_{FX}Z,Y)\\
			&=g(h(Y,Z),FX)
			\intertext{which gives \eqref{eq:g(h(D_1,D_2),FD_1) Symmetric Bi-Slant Warped Product Submanifolds}. Repeating the above calculation, we have}
			g(h(X,Z),FY)&=g(\overline{\nabla}_XZ,JY-PY)\\
			&=-g(J\overline{\nabla}_XZ,Y)-g(\nabla_XZ,PY)\\
			&=-g(\overline{\nabla}_XJZ,Y)-\frac{1}{2}g(JB,Z)g(X,Y)-\frac{1}{2}g(B,Z)g(JX,Y)\\
			&=-g(\overline{\nabla}_XPZ,Y)-g(\overline{\nabla}_XFZ,Y)-\frac{1}{2}g(JB,Z)g(X,Y)\\
			&\hspace{0.4cm}-\frac{1}{2}g(B,Z)g(JX,Y)\\
			&=g(\mathfrak{A}_{FZ}X,Y)+\frac{1}{2}g(B,JZ)g(X,Y)-\frac{1}{2}g(B,Z)g(PX,Y)\\
			\intertext{Using \eqref{eq:g(h(D_1,D_2),FD_1) Symmetric Bi-Slant Warped Product Submanifolds} and comparing symmetric and skew symmetric terms in $X$ and $Y$ we have,}
			\frac{1}{2}g(B,Z)&=0
			\intertext{which gives \eqref{eq:Identity g(B,D_2) Bi-Slant Warped Product Submanifolds} and }
			g(h(X,Z),FY)&=g(h(X,Y),FZ)+\frac{1}{2}g(B,PZ+FZ)g(X,Y)
		\end{align*}
		which on substituting from \eqref{eq:Identity g(B,D_2) Bi-Slant Warped Product Submanifolds} gives \eqref{eq:Identity g(h(D_1,D_1),FD_2) & g(h(D_1,D_2),FD_1) Bi-Slant Warped Product Submanifolds}.
	\end{proof}
	\noindent From \eqref{eq:Identity g(B,D_1) Bi-Slant Warped Product Submanifolds} and \eqref{eq:Identity g(B,D_2) Bi-Slant Warped Product Submanifolds} we have 
	\begin{cor}
		Given a pointwise bi-slant warped product submanifold \mbox{$M= M_{\theta_1}\times \,_{\lambda}M_{\theta_2}$} in an lcK manifold $(\widetilde{M}^{2n},J,g )$, we have the Lee vector field $B$ is orthogonal to the second factor and the warping function $\lambda$ satisfies $\text{grad}(\ln\lambda)=\frac{1}{2}B^T$, where $B^T$ denotes the tangential part of the Lee vector field along $M$.
	\end{cor}
	\begin{rem}\label{rem:Local Orthonormal Basis of Bi-Slant Warped Product Submanifolds}
		Given a pointwise bi-slant warped product submanifold $M_{\theta_1}\times\,_\lambda M_{\theta_2}$ of an l.c.K manifold $\widetilde{M}^{2n}$, let $\{X_i,\beta_1X_i\}_{i=1}^{p}$ and $\{Z_j,\beta_2PZ_j\}_{j=1}^{q}$ respectively be local orthonormal frames of $TM_{\theta_1}$ and $TM_{\theta_2}$. Then a local orthonormal frame of $\widetilde{M}^{2n}$ is
		\[\begin{aligned}
			&\left\lbrace \widehat{X_i}=X_i,\widehat{PX_i}=\beta_1PX_i\right\rbrace\cup\left\lbrace\widehat{Z_j}=\frac{Z_j}{\lambda},\widehat{PZ_j}=\frac{\beta_2 PZ_j}{\lambda}\right\rbrace\cup\left\lbrace\widehat{\xi_k},\widehat{J\xi_k}\right\rbrace\\&\cup\left\lbrace\widehat{FX_i}=\alpha_1 FX_i,\widehat{FPX_i}=\alpha_1\beta_1 FPX_i\right\rbrace\cup\left\lbrace\widehat{FZ_j}=\frac{\alpha_2 FZ_j}{\lambda},\widehat{FPZ_j}=\frac{\alpha_2\beta_2 FPZ_j}{\lambda}\right\rbrace
		\end{aligned}\]
		where $\alpha_i=\csc\theta_i$, $\beta_i=\sec\theta_i$ for $i=1,2$ and
		\begin{align*}
			&\left\{\widehat{X_i}\widehat{PX_i}:1\leq i\leq n_1\right\}\text{ is an orthonormal basis of }\mathcal{D}^{\theta_1}\\
			&\left\{\widehat{Z_j}, \widehat{PZ_j}:1\leq j\leq n_2\right\}\text{ is an orthonormal basis of }\mathcal{D}^{\theta_2}\\
			&\left\{\widehat{FX_i}, \widehat{FPX_i}:1\leq j\leq n_1\right\}\text{ is an orthonormal basis of }F\mathcal{D}^{\theta_1}\\
			&\left\{\widehat{FZ_j}, \widehat{FPZ_j}:1\leq j\leq n_2\right\}\text{ is an orthonormal basis of }F\mathcal{D}^{\theta_2}\\
			&\left\{\widehat{\xi_k},\widehat{J\xi_k}:1\leq s\leq \frac{n-2n_1-2n_2}{2}\right\}\text{ is an orthonormal basis of }\mu
		\end{align*}
		However, while $Z_j,\beta PZ_j\in\mathcal{L}(M_{\theta_2})$ we have $\widehat{Z_j},\widehat{PZ_j}\notin\mathcal{L}(M_{\theta_2})$ in general, as $\lambda$ is a function on $M_{\theta_1}$. Also, note that
		\begin{align*}
			J\left(\widehat{Z_j}\right)&=J\left(\frac{Z_j}{\lambda}\right)=\frac{PZ_j}{\lambda}+\frac{FZ_j}{\lambda}=\cos\theta_2\widehat{PZ_j}+\sin\theta_2\widehat{FZ_j}\\
			J\left(\widehat{PZ_j}\right)&=J\left(\sec\theta_2\frac{PZ_j}{\lambda}\right)=\frac{\sec\theta_2 P^2Z_j}{\lambda}+\frac{\sec\theta_2 FPZ_j}{\lambda}=-\cos\theta_2\widehat{Z_j}+\sin\theta_2\widehat{FPZ_j}
		\end{align*}
	\end{rem}
	\section{Main Results}
	We first give a characterisation for pointwise bi-slant warped product submanifolds of lcK manifolds.
	\begin{theo}\label{th:Bi-Slant Warped Prod Characterisation}
		Let $M^m$ be a pointwise bi-slant submanifold of an lcK manifold $\widetilde{M}^{2n}$. Then the following are equivalent 
		\begin{enumerate}
			\item $M^m$ is a pointwise bi-slant warped product submanifold \mbox{$M_{\theta_1}\times\,_\lambda M_{\theta_2}$} of $\widetilde{M}^{2n}$\label{Bi-Slant Warped Prod Characterisation}
			\item For every $X,Y\in\mathcal{L}(M_{\theta_1})$ and $Z,W\in\mathcal{L}(M_{\theta_2})$ we have
			\begin{align}\label{eq:Bi-Slant Warped Prod Characterisation 1}
				\omega(\mathcal{D}^{\theta_2})=\{0\}\text{ and }g\left(\mathfrak{A}_{FPX}Z\right.&\left.-\ \mathfrak{A}_{FX}PZ,Y\right)+g(\nabla^\perp_YFX,FZ)=0\nonumber\\
					g(\mathfrak{A}_{FPZ}X-\mathfrak{A}_{FZ}PX,W)&+g(\nabla^\perp_WFZ,FX)\\
					&=\sin^2\theta_2\left(\frac{1}{2}\omega(X)-X(\ln\lambda)\right)g(Z,W)\nonumber
			\end{align}
			for some smooth function $\lambda:M_{\theta_1}\to(0,\infty)$.\label{Bi-Slant Warped Prod Characterisation - Condition 1}
			\item For every $X\in\mathcal{L}(M_{\theta_1})$ and $Z\in\mathcal{L}(M_{\theta_2})$ we have
			\begin{equation}\label{eq:Bi-Slant Warped Prod Characterisation 2}
				\omega(\mathcal{D}^{\theta_2})=\{0\}\text{ and }\nabla_XZ=\nabla_ZX=\frac{1}{2}\omega(X)Z
			\end{equation}\label{Bi-Slant Warped Prod Characterisation - Condition 2}
		\end{enumerate}
		Also, in this case we have the mean curvature vector $\mathcal{H}$ of $M_{\theta_2}$ in $M^m$ is 
		\begin{equation}\label{eq:Warping Function Bi-Slant Warped Prod Characterisation}
			\mathcal{H}=-\text{grad}(\ln\lambda)=-\frac{1}{2}B^T
		\end{equation}
		where $B^T$ is the tangential component of $B$ along $M$.
	\end{theo}
	\begin{proof}
		\mbox{\eqref{Bi-Slant Warped Prod Characterisation}$\Leftrightarrow$\eqref{Bi-Slant Warped Prod Characterisation - Condition 1}} This follows from Theorem \ref{th:Bi-Slant D_theta1}, Theorem \ref{th:Bi-Slant D_theta2} and $\text{grad}(\ln\lambda)\in\mathcal{L}(M_{\theta_1})$ which implies
		\begin{align*}
			g(\nabla_Z(\text{grad}(\ln\lambda)),X)&=ZX(\ln\lambda)-g(\text{grad}(\ln\lambda),\nabla_ZX)\\
			&=[Z,X](\ln\lambda)-\nabla_ZX(\ln\lambda)\\
			&=-\nabla_XZ(\ln\lambda)\\
			&=g(Z,\nabla_X(\text{grad}(\ln\lambda)))\\
			&=0\hspace{0.4cm}
		\end{align*}
		as $Z(\ln\lambda)=0$ and $M_{\theta_1}$ is totally geodesic in $M$. Also, \eqref{eq:Warping Function Bi-Slant Warped Prod Characterisation} follows from Lemma \ref{lem:Identities For Bi-Slant Warped Product Submanifolds} \eqref{eq:Identity g(B,D_1) Bi-Slant Warped Product Submanifolds}.\par\noindent
		\mbox{\eqref{Bi-Slant Warped Prod Characterisation}$\Leftrightarrow$\eqref{Bi-Slant Warped Prod Characterisation - Condition 2}} Let \mbox{$M=M_{\theta_1}\times\,_\lambda M_{\theta_2}$} be a pointwise bi-slant warped product submanifold. Then \eqref{eq:Bi-Slant Warped Prod Characterisation 2} and \eqref{eq:Warping Function Bi-Slant Warped Prod Characterisation} follow from \eqref{eq:Connection Warped Product on M_1 x M_2} and Lemma \ref{lem:Identities For Bi-Slant Warped Product Submanifolds} \eqref{eq:Identity g(B,D_1) Bi-Slant Warped Product Submanifolds}.\par
		Conversely, let $M^m$ be a pointwise bi-slant submanifold of an lcK manifold $\widetilde{M}^{2n}$ such that \eqref{eq:Bi-Slant Warped Prod Characterisation 2} holds. Then for all $X,Y\in \mathcal{L}(M_{\theta_1})$ and $Z,W\in \mathcal{L}(M_{\theta_2})$ we have
		\begin{align*}
			g([X,Y],Z)&=g(\nabla_XY-\nabla_YX,Z)\\
			&=-g(\nabla_XZ,Y)+g(\nabla_YZ,X)\\
			&=0
			\intertext{which implies $\mathcal{D}^{\theta_1}$ is involutive.}
			g(\nabla_XY,Z)&=-g(\nabla_XZ,Y)=0
			\intertext{which implies leaves of $\mathcal{D}^{\theta_1}$ are totally geodesic in $M$.}
			g([Z,W],X)&=g(\nabla_ZW-\nabla_WZ,X)\\
			&=-g(\nabla_ZX,W)+g(\nabla_WX,Z)\\
			&=-\frac{1}{2}\omega(X)g(Z,W)+\frac{1}{2}\omega(X)g(W,Z)=0
			\intertext{which implies $\mathcal{D}^{\theta_2}$ is involutive.}
			g(\nabla_ZW,X)&=-g(\nabla_ZX,W)\\
			&=-\frac{1}{2}\omega(X)g(Z,W)\\
			&=-\frac{1}{2}g(Z,W)g(B^T,X)
			\intertext{which implies leaves of $\mathcal{D}^{\theta_2}$ are totally umbilical in $M$ with mean curvature vector $-\frac{1}{2}B^T$.}
			g\left(\nabla_ZB^T,X\right)&=\frac{1}{2}\omega\left(B^T\right)g(Z,X)=0
		\end{align*}
		which implies $B^T$ is parallel in the normal bundle of $M_{\theta_2}$ in $M$.\par
		Hence by Theorem \ref{th:Characterisation - Warped Product} we have \mbox{$M=M_{\theta_1}\times\,_\lambda M_{\theta_2}$} is a pointwise bi-slant warped product submanifold.
	\end{proof}
	We conclude our study of pointwise bi-slant warped product submanifolds of lcK manifolds by giving an inequality for the squared norm of the second fundamental form.
	\begin{theo}\label{th:Length h - Bi-Slant Warped Product}
		Let \mbox{$M= M_{\theta_1}\times \,_\lambda M_{\theta_2}$} be a pointwise bi-slant warped product submanifold in an lcK manifold $(\widetilde{M}^{2n},J,g )$. Then the norm of the second fundamental form satisfies the inequality
		\begin{align}\label{eq:Length h - Bi-Slant Warped Product}
			||h||^2\geq\:&
			\frac{n_1}{2}\sin^2\theta_2\|B|_{F\mathcal{D}^{\theta_2}}\|^2+\frac{n_2}{2}\sin^2\theta_1\|B|_{F\mathcal{D}^{\theta_1}}\|^2+\sin\theta_2g\left(H_{\mathcal{D}^{\theta_1}}|_{F\mathcal{D}^{\theta_2}},B|_{F\mathcal{D}^{\theta_2}}\right)\nonumber\\&\hspace{0.4cm}+\sin\theta_1g\left(H_{\mathcal{D}^{\theta_2}}|_{F\mathcal{D}^{\theta_1}},B|_{F\mathcal{D}^{\theta_1}}\right)
		\end{align}
		where $2n_1=\dim_{\mathbb{R}}\mathcal{D}^{\theta_1}$, $2n_2=\dim_{\mathbb{R}}\mathcal{D}^{\theta_2}$ and $H_{\mathcal{D}^{\theta_1}}$ and $H_{\mathcal{D}^{\theta_2}}$ are respectively the components of the mean curvature vector $H$ of $M$ in $\widetilde{M}^{2n}$ along $\mathcal{D}^{\theta_1}$ and $\mathcal{D}^{\theta_2}$.\par
		\noindent If equality holds then we have
		\begin{itemize}
			\item $\text{Image}(h)\subseteq(F\mathcal{D}^{\theta_1}\oplus F\mathcal{D}^{\theta_2})$, and
			\item $M$ is minimal in $\widetilde{M}^{2n}$, if and only if, $M$ is mixed-totally geodesic in $\widetilde{M}^{2n}$.
		\end{itemize}
	\end{theo}
	\begin{proof}
		\begin{align*}
			||h||^2=&\left\|h(\mathcal{D}^{\theta_1},\mathcal{D}^{\theta_1})\big|_{F\mathcal{D}^{\theta_1}}\right\|^2
			+\left\|h(\mathcal{D}^{\theta_1},\mathcal{D}^{\theta_2})\big|_{F\mathcal{D}^{\theta_1}}\right\|^2
			+\left\|h(\mathcal{D}^{\theta_2},\mathcal{D}^{\theta_2})\big|_{F\mathcal{D}^{\theta_1}}\right\|^2\\
			&+\left\|h(\mathcal{D}^{\theta_1},\mathcal{D}^{\theta_1})\big|_{F\mathcal{D}^{\theta_2}}\right\|^2
			+\left\|h(\mathcal{D}^{\theta_1},\mathcal{D}^{\theta_2})\big|_{F\mathcal{D}^{\theta_2}}\right\|^2
			+\left\|h(\mathcal{D}^{\theta_2},\mathcal{D}^{\theta_2})\big|_{F\mathcal{D}^{\theta_2}}\right\|^2\\
			&+\left\|h(\mathcal{D}^{\theta_1},\mathcal{D}^{\theta_1})\big|_{\mu}\right\|^2
			+\left\|h(\mathcal{D}^{\theta_1},\mathcal{D}^{\theta_2})\big|_{\mu}\right\|^2
			+\left\|h(\mathcal{D}^{\theta_2},\mathcal{D}^{\theta_2})\big|_{\mu}\right\|^2
		\end{align*}
		From \eqref{eq:Identity g(h(D_1,D_1),FD_2) & g(h(D_1,D_2),FD_1) Bi-Slant Warped Product Submanifolds} and Remark \ref{rem:Local Orthonormal Basis of Bi-Slant Warped Product Submanifolds} we have
		\begin{align*}
			g\left( h\left(\widehat{X_i},\widehat{Z_p}\right),\widehat{FX_j}\right)
			=&\frac{\csc\theta_1}{\lambda}\left\lbrace
			g\left( h\left(X_i,X_j\right),FZ_p\right)+\frac{1}{2}\delta_{ij}g\left( B,FZ_p\right)
			\right\rbrace\\
			\implies\sin\theta_1g\left( h\left(\widehat{X_i},\widehat{Z_p}\right),\widehat{FX_j}\right)
			=&\sin\theta_2g\left( h\left(\widehat{X_i},\widehat{X_j}\right),\widehat{FZ_p}\right)
			+\frac{1}{2}\sin\theta_2\delta_{ij}g\left( B,\widehat{FZ_p}\right)
		\end{align*}
		Similarly,
		\begin{align*}
			\sin\theta_1g\left( h\left(\widehat{X_i},\widehat{PZ_p}\right),\widehat{FX_j}\right)
			=&\sin\theta_2g\left( h\left(\widehat{X_i},\widehat{X_j}\right),\widehat{FPZ_p}\right)
			+\frac{1}{2}\sin\theta_2\delta_{ij}g\left( B,\widehat{FPZ_p}\right)\\
			\sin\theta_1g\left( h\left(\widehat{PX_i},\widehat{Z_p}\right),\widehat{FPX_j}\right)
			=&\sin\theta_2g\left( h\left(\widehat{PX_i},\widehat{PX_j}\right),\widehat{FZ_p}\right)
			+\frac{1}{2}\sin\theta_2\delta_{ij}g\left( B,\widehat{FZ_p}\right)\\
			\sin\theta_1g\left( h\left(\widehat{PX_i},\widehat{PZ_p}\right),\widehat{FPX_j}\right)
			=&\sin\theta_2g\left( h\left(\widehat{PX_i},\widehat{PX_j}\right),\widehat{FPZ_p}\right)\\
			&+\frac{1}{2}\sin\theta_2\delta_{ij}g\left( B,\widehat{FPZ_p}\right)\\
			\sin\theta_1g\left( h\left(\widehat{PX_i},\widehat{Z_p}\right),\widehat{FX_j}\right)
			=&\sin\theta_2g\left( h\left(\widehat{PX_i},\widehat{X_j}\right),\widehat{FZ_p}\right)\\
			\sin\theta_1g\left( h\left(\widehat{PX_i},\widehat{PZ_p}\right),\widehat{FX_j}\right)
			=&\sin\theta_2g\left( h\left(\widehat{PX_i},\widehat{X_j}\right),\widehat{FPZ_p}\right)\\
			\sin\theta_1g\left( h\left(\widehat{X_i},\widehat{Z_p}\right),\widehat{FPX_j}\right)
			=&\sin\theta_2g\left( h\left(\widehat{X_i},\widehat{PX_j}\right),\widehat{FZ_p}\right)\\
			\sin\theta_1g\left( h\left(\widehat{X_i},\widehat{PZ_p}\right),\widehat{FPX_j}\right)
			=&\sin\theta_2g\left( h\left(\widehat{X_i},\widehat{PX_j}\right),\widehat{FPZ_p}\right)
		\end{align*}
		which implies
		\begin{align*}
			\left\|h(\mathcal{D}^{\theta_1},\mathcal{D}^{\theta_2})\big|_{F\mathcal{D}^{\theta_1}}\right\|^2
			=&\cos^2\theta_1\left\|h(\mathcal{D}^{\theta_1},\mathcal{D}^{\theta_2})\big|_{F\mathcal{D}^{\theta_1}}\right\|^2
			+\sin^2\theta_1\left\|h(\mathcal{D}^{\theta_1},\mathcal{D}^{\theta_2})\big|_{F\mathcal{D}^{\theta_1}}\right\|^2\\
			=&\cos^2\theta_1\left\|h(\mathcal{D}^{\theta_1},\mathcal{D}^{\theta_2})\big|_{F\mathcal{D}^{\theta_1}}\right\|^2
			+\sin^2\theta_2\left\|h(\mathcal{D}^{\theta_1},\mathcal{D}^{\theta_1})\big|_{F\mathcal{D}^{\theta_2}}\right\|^2\\
			&+\frac{1}{2}\sin^2\theta_2\sum_{i,p}\left\lbrace 
			g\left( B,\widehat{FZ_p}\right)^2
			+g\left( B,\widehat{FPZ_p}\right)^2
			\right\rbrace\\
			&+\sin\theta_2\sum_{i,p}\left\lbrace
			g\left( h\left(\widehat{X_i},\widehat{X_i}\right),\widehat{FZ_p}\right)g\left( B,\widehat{FZ_p}\right)\right.\\
			&\hspace{2.4cm}\left.+g\left( h\left(\widehat{X_i},\widehat{X_i}\right),\widehat{FPZ_p}\right)g\left( B,\widehat{FPZ_p}\right)\right.\\
			&\hspace{2.4cm}\left.+g\left( h\left(\widehat{PX_i},\widehat{PX_i}\right),\widehat{FZ_p}\right)g\left( B,\widehat{FZ_p}\right)\right.\\
			&\hspace{2.4cm}\left.+g\left( h\left(\widehat{PX_i},\widehat{PX_i}\right),\widehat{FPZ_p}\right)g\left( B,\widehat{FPZ_p}\right)
			\right\rbrace
		\end{align*}
		i.e.
		\begin{align*}
			\sin^2\theta_1\left\|h(\mathcal{D}^{\theta_1},\mathcal{D}^{\theta_2})\big|_{F\mathcal{D}^{\theta_1}}\right\|^2
			&=\sin^2\theta_2\left\|h(\mathcal{D}^{\theta_1},\mathcal{D}^{\theta_1})\big|_{F\mathcal{D}^{\theta_2}}\right\|^2
			+\frac{2n_1}{4}\sin^2\theta_2\biggl\|B\big|_{F\mathcal{D}^{\theta_2}}\biggr\|^2\\
			&+\sin\theta_2g\left(\left.\sum_i\left\lbrace
			h\left(\widehat{X_i},\widehat{X_i}\right)+h\left(\widehat{PX_i},\widehat{PX_i}\right)
			\right\rbrace\right|_{F\mathcal{D}^{\theta_2}},B|_{F\mathcal{D}^{\theta_2}}\right)\\
			&\geq\frac{n_1}{2}\sin^2\theta_2\|B|_{F\mathcal{D}^{\theta_2}}\|^2+\sin\theta_2g\left(H_{\mathcal{D}^{\theta_1}}|_{F\mathcal{D}^{\theta_2}},B|_{F\mathcal{D}^{\theta_2}}\right)
		\end{align*}
		As done above we have
		\begin{align*}
			\sin^2\theta_2\left\|h(\mathcal{D}^{\theta_1},\mathcal{D}^{\theta_2})\big|_{F\mathcal{D}^{\theta_2}}\right\|^2
			&=\sin^2\theta_1\left\|h(\mathcal{D}^{\theta_2},\mathcal{D}^{\theta_2})\big|_{F\mathcal{D}^{\theta_1}}\right\|^2
			+\frac{2n_2}{4}\sin^2\theta_1\biggl\|B\big|_{F\mathcal{D}^{\theta_1}}\biggr\|^2\\
			&+\sin\theta_1g\left(\left.\sum_p\left\lbrace
			h\left(\widehat{Z_p},\widehat{Z_p}\right)+h\left(\widehat{PZ_p},\widehat{PZ_p}\right)
			\right\rbrace\right|_{F\mathcal{D}^{\theta_1}},B|_{F\mathcal{D}^{\theta_1}}\right)\\
			&\geq\frac{n_2}{2}\sin^2\theta_1\|B|_{F\mathcal{D}^{\theta_1}}\|^2+\sin\theta_1g\left(H_{\mathcal{D}^{\theta_2}}|_{F\mathcal{D}^{\theta_1}},B|_{F\mathcal{D}^{\theta_1}}\right)
		\end{align*}
		Combining we have \eqref{eq:Length h - Bi-Slant Warped Product}.\par
		If equality holds in \eqref{eq:Length h - Bi-Slant Warped Product}, then the only non-zero components of $||h||$ are
		\begin{align*}
			\left\|h(\mathcal{D}^{\theta_1},\mathcal{D}^{\theta_1})|_{F\mathcal{D}^{\theta_2}}\right\|^2\hspace{0.25cm}\text{  , }\hspace{0.35	cm} \left\|h(\mathcal{D}^{\theta_1},\mathcal{D}^{\theta_2})|_{F\mathcal{D}^{\theta_1}}\right\|^2,\\ \left\|h(\mathcal{D}^{\theta_2},\mathcal{D}^{\theta_2})|_{F\mathcal{D}^{\theta_1}}\right\|^2\hspace{0.1cm}\text{and }\hspace{0.1cm} \left\|h(\mathcal{D}^{\theta_1},\mathcal{D}^{\theta_2})|_{F\mathcal{D}^{\theta_2}}\right\|^2.
		\end{align*}
		Also, from the calculations above we have, 
		\begin{align*}
			\left\|h(\mathcal{D}^{\theta_1},\mathcal{D}^{\theta_1})|_{F\mathcal{D}^{\theta_2}}\right\|^2=0&\text{ if and only if }\left\|h(\mathcal{D}^{\theta_1},\mathcal{D}^{\theta_2})|_{F\mathcal{D}^{\theta_1}}\right\|^2=0\\&\hspace{0.5cm}\text{and}\\
			\left\|h(\mathcal{D}^{\theta_2},\mathcal{D}^{\theta_2})|_{F\mathcal{D}^{\theta_1}}\right\|^2=0&\text{ if and only if }\left\|h(\mathcal{D}^{\theta_1},\mathcal{D}^{\theta_2})|_{F\mathcal{D}^{\theta_2}}\right\|^2=0 
		\end{align*}
		Hence, the result follows.
	\end{proof}
	\section{Example}
	Consider $\mathbb{C}^n=\mathbb{E}^{2n},J,g_0)$ where $\mathbb{E}^{2n}$ is the Euclidean space of dimension $2n$ with coordinates $(x_1,\ldots,x_n,y_1,\ldots,y_n)$ equipped with the standard Euclidean metric $g_0$ and the canonical almost complex structure
	\begin{equation}\label{Example:J}
		J(x_1,\ldots,x_n,y_1,\ldots,y_n)=(-y_1,\ldots,-y_n,x_1,\ldots,x_n)
	\end{equation}
	Then $\mathbb{C}^n=(\mathbb{E}^{2n},J,g_0)$ is a flat Kähler manifold.\par
	We use the following result about pointwise slant immersions.
	\begin{theo}[\cite{Chen-PointwiseSlant} (Proposition 2.2)]\label{Example:theorem}
		Given a Kähler manifold $(\widetilde{M}^{2n},J,g_0)$, let $M_\theta$ be a pointwise slant submanifold. Then for any smooth function $f:\widetilde{M}\to(0,\infty)$, we have $M_\theta$ is again a pointwise slant submanifold of the globally conformal Kähler (gcK) manifold $(\widetilde{M}^{2n},J,e^{-f}g_0)$ with the same slant angle.
	\end{theo}
	\begin{eg}
		Let $\mathbb{C}^4=(\mathbb{E}^8,J,g_0)$ be as defined above. Consider an open subset of $\mathbb{E}^4$ with $u_1u_2\neq1$, $u_3u_4\neq1$, $(u_1-u_2)\in\left(0,\frac{\pi}{4}\right)$ and  $(u_3-u_4)\in\left(\frac{\pi}{4},\frac{\pi}{2}\right)$. Define the $4$-dimensional submanifold $M$ of $\mathbb{C}^4$ given by
		\begin{equation}\label{Example:bislant}
			\begin{aligned}
				x_1&=u_1\cos u_2&,&&y_1&=u_1\sin u_2\\
				x_2&=u_2\cos u_1&,&&y_2&=u_2\sin u_1\\
				x_3&=u_3\cos u_4&,&&y_3&=u_3\sin u_4\\
				x_4&=u_4\cos u_3&,&&y_4&=u_4\sin u_3
			\end{aligned}
		\end{equation}
		An orthonormal frame of the tangent bundle $TM$ of $M$ is
		\begin{align*}
			X_1&=\frac{1}{\sqrt{1+u_2^2}}\left(\cos u_2\frac{\partial}{\partial x_1}-u_2\sin u_1\frac{\partial}{\partial x_2}+\sin u_2\frac{\partial}{\partial y_1}+u_2\cos u_1\frac{\partial}{\partial y_2}\right)\\
			X_2&=\frac{1}{\sqrt{1+u_1^2}}\left(-u_1\sin u_2\frac{\partial}{\partial x_1}+\cos u_1\frac{\partial}{\partial x_2}+u_1\cos u_2\frac{\partial}{\partial y_1}+\sin u_1\frac{\partial}{\partial y_2}\right)\\
			X_3&=\frac{1}{\sqrt{1+u_4^2}}\left(\cos u_4\frac{\partial}{\partial x_3}-u_4\sin u_3\frac{\partial}{\partial x_4}+\sin u_4\frac{\partial}{\partial y_3}+u_4\cos u_3\frac{\partial}{\partial y_4}\right)\\
			X_4&=\frac{1}{\sqrt{1+u_3^2}}\left(-u_3\sin u_4\frac{\partial}{\partial x_3}+\cos u_3\frac{\partial}{\partial x_4}+u_3\cos u_4\frac{\partial}{\partial y_3}+\sin u_3\frac{\partial}{\partial y_4}\right)
		\end{align*}
		Then, $M$ is a proper pointwise bi-slant submanifold with slant distributions given by $\mathcal{D}^{\theta_1}=\text{Span}\{X_1,X_2\}$ and $\mathcal{D}^{\theta_2}=\text{Span}\{X_3,X_4\}$. Also, the slant angles are given by
		\[\cos^2\theta_1=\frac{(u_1u_2-1)^2\cos^2(u_1-u_2)}{(1+u_1^2)(1+u_2^2)}\text{ , and , }\cos^2\theta_2=\frac{(u_3u_4-1)^2\cos^2(u_3-u_4)}{(1+u_3^2)(1+u_4^2)}\]
		It is straightforward to check that $\mathcal{D}^{\theta_1}$ and $\mathcal{D}^{\theta_2}$ are both involutive and totally geodesic in $M$. Let $M_{\theta_1}$ and $M_{\theta_2}$ be the leaves of $\mathcal{D}^{\theta_1}$ and $\mathcal{D}^{\theta_2}$ respectively. Then we have $M$ is the Riemannian product $M=M_{\theta_1}\times M_{\theta_2}$ and the metric $g_M$ induced on $M$ from $\mathbb{C}^4$ is given by
		\begin{equation}\label{Example:Induced metric M}
			g_M=g_1+g_2
		\end{equation}
		where
		\begin{equation}\label{Example:Induced metrics M_1 and M_2}
			g_1=(1+u_2^2)du_1^2+(1+u_1^2)du_2^2\text{ , and , }g_2=(1+u_4^2)du_3^2+(1+u_3^2)du_4^2
		\end{equation}
		Now, for any non-constant positive smooth function $f=f(x_1,x_2,y_1,y_2)$ on $\mathbb{C}^4$, depending only on coordinates $x_1,x_2,y_1,y_2$, consider the Riemannian metric $\tilde{g}=e^{-f}g_0$, conformal to the standard metric $g_0$. Then, $\widetilde{M}=(\mathbb{E}^8,J,\tilde{g})$ is a globally conformal Kähler manifold and the metric on $M$ induced from $\widetilde{M}$ is the warped product metric
		\begin{equation}\label{Example:Induced warped metric M}
			\tilde{g}_M=\tilde{g}_1+e^{-f}g_2
		\end{equation}
		where
		\begin{equation}\label{Example:Induced warped metrics M_1 and M_2}
			\tilde{g}_1=e^{-f}g_1
		\end{equation}
		is conformal to $g_1$ by the choice of $f$.\par
		Hence from Theorem \ref{Example:theorem} we have $(M,\tilde{g}_M)$ is a proper pointwise bi-slant warped product submanifold of $\widetilde{M}=(\mathbb{E}^8,J,\tilde{g})$.\par
		Also, as  $f=f(x_1,x_2,y_1,y_2)$ is a non-constant positive smooth function on $\mathbb{C}^4$, depending only on coordinates $x_1,x_2,y_1,y_2$, from \eqref{Example:bislant} we have that restricted to the submanifold $M$, the Lee form $\omega$ of $\widetilde{M}$ is given by
		\begin{equation}\label{Example:B}
			\omega=df=\frac{\partial f}{\partial u_1}du_1+\frac{\partial f}{\partial u_2}du_2
		\end{equation}
		Hence, it follows that the Lee-vector field $B$ is orthogonal to $\mathcal{D}^{\theta_2}$ and the warping function $\lambda=-e^{-\frac{f}{2}}|_M$ satisfies $\text{grad}(\ln \lambda)=\frac{1}{2}\text{grad}(f|_M)=\frac{1}{2}B^T$.
	\end{eg}
	
\end{document}